\documentstyle[11pt]{article}
\input amssym.def
\input amssym
\def\bpi{\pi_{\R^2}}
 \def\d{ \rm d}

\def\l{\lambda}

\def\tw{\rm tw}

\def\w{\wedge}

\def\L{\Lambda}
\def\G{\Gamma}
\def\W{\bigwedge}
\newcommand{\R}{\mbox{${\Bbb R}$}}

\newcommand{\C}{\mbox{${\Bbb C}$}}

\newcommand{\F}{I\!\!F}
\def\CN{\cal C}

\def\Ham{\rm Ham}

\def\Diff{\rm Diff }
\def\max{\scriptscriptstyle {\rm max} }
\def\ie{i.e.\ }

\def\si{\sigma}
\def\pa{\partial}
\def\lw{\longrightarrow}
\def\lo{\longmapsto}

\newcommand{\fg}{\mbox{${\frak g}$}}
\newcommand{\fh}{\mbox{${\frak h}$}}

\def\cf{\mbox{$\rm cf$}}

\def\la{\langle}
\def\ra{\rangle}

\def\dim{\mbox{$\rm dim$}}

\newcommand{\fl}{\mbox{${\cal F}$}}
\newcommand{\sh}{\mbox{${\sharp}$}}
\def\sect#1{\section*{\centerline{\large\bf  #1}}}
\def\jielun#1{{\par\noindent\hskip 0.5 cm}
{\bf\heiti #1}\bgroup\it}
\def\endjielun{\egroup\par\bigskip}
\def\bli#1{{\par\noindent\hskip 0.5 cm}
{\bf #1}\bgroup\rm}
\def\eli{\egroup\par\bigskip}
\newcommand{\qed}{\hfill $\Box$ }
\newcommand {\pf}{{\par\noindent\hskip 0.65 cm}{\it  Proof. }}
\begin{document}
\title{\bf   {\CC {\char64}n K+ 4z
J{\char125} Sk A? WS 4z J{\char125} 5D JU Ku}}
\title{{\bf {Almost Homogeneous Poisson Spaces}}
\thanks{Research partially supported by  SFB-237 of the DFG.}
}
\author{Qi-Lin Yang\\
Department of Mathematics,
Tsinghua University, \\
Beijing 100084, P. R. China\\
(E-mail: qlyang@math.tsinghua.edu.cn)\\
  }
\maketitle
\date{}

    \begin{abstract}
We prove that any holomorphic Poisson manifold has an open
symplectic leaf which is a  pseudo-K\"ahler submanifold, and we
define an  obstruction to study the equivariance of momentum map
for tangential Poisson action.  Some properties of almost
homogeneous Poisson manifolds are studied and we show that any
compact symplectic Poisson homogeneous space is a torus bundle
over a dressing orbit.
\end{abstract}
\sect{\S1. \  Introduction}
\hskip\parindent
Motivated by the study of Poisson nature of ``dressing
transformation"  for soliton equations, Lu  and Weinstein
  introduced the notion of Poisson action {[STS, Lu, LW, W]}.
  Poisson actions  are used  to understand the
   ``hidden symmetries"  of certain integrable systems.
When the Poisson structure of the Poisson Lie group is  zero, the
Poisson structure of the acted Poisson manifold is preserved. So
Poisson actions are the
 generalizations of symplectic actions and Hamiltonian
 actions. Tangential
  Poisson action $G\times P \lw P$ is a  special
 kind of Poisson action which is the most adjacent generalization
   of symplectic action. The orbits of tangential
 Poisson action  lie in the symplectic leaves of $P$
and the restricted
    action of $G$ to every symplectic
 leaf is also a Poisson action  (\cf. [Y]).

     Momentum mapping for symplectic  action on
 symplectic manifold is the abstract of momentum and angular
   momentum in classical mechanics.The symmetry of
 the phase space of a mechanical system leads to
 some conservative quantities,
   to reduce the symmetry by these  conservative
 quantities is a way to simplify the mechanical system.
    So the momentum map plays an important role in
reduction. In {[Lu]}, Lu  generalized  this concept to the
tangential Poisson
    action and studied the Meyer-Marsden-Weinstein
momentum  reduction for tangential Poisson action on
 symplectic manifold in regular case. In
    [Y], we generalized momentum reduction to singular
value on any Poisson manifold.
    In [Gi], Ginzburg studied  the further properties
of momentum map for tangential Poisson action.

   In this paper, we discuss  other properties of
the momentum  mapping for tangential Poisson actions.
In sect. 2, we give a detailed description of symplectic
stratification of a Poisson manifold.
    In sect. 3,  we  in principle classify  Poisson
Lie structures on abelian Lie groups.
  In sect. 4, we give the  necessary and sufficient
condition when a Poisson action is Poisson structure
 preserved, and define an obstruction to describe
the equivariance of momentum map.
      In sect. 5,  we prove that the  compact symplectic
Poisson homogeneous space is a torus bundle over some dressing
orbit, and give an example  of non-transitive tangential Poisson
action which is not a dressing action.

\sect{\S2. \  Poisson manifolds and their symplectic
stratifications} \hskip 0.5cm A commutative associative algebra
${\cal A}$ over field $\F$ is called a {\it Poisson algebra} if it
is equipped with a $\F$-bilinear operation (called {\it Poisson
bracket})
 $\{,\}: {\cal A}\times {\cal A}\lw {\cal A}$
satisfying the following conditions:

{\rm (i)} $({\cal A}, \{,\})$ is a Lie algebra over $\F;$

{\rm (ii)} the commutative algebraic structure is compatible with
the Lie algebraic structure in the following way, $\ie,$ the
Leibniz rule is satisfied:
$$\{ab,c\}=a\{b,c\}+\{a,c\}b.$$

  An algebraic $\F$-variety $V$ is called
an {\it algebraic Poisson $\F$-variety} if the sheaf ${\cal
O}_{V}$ of regular functions carries  a structure of a sheaf of
Poisson algebra. A complex space $W$ is called a Poisson space if
the sheaf of regular functions is a holomorphic Poisson algebra,
in the smooth case $W$ is called a {\it  holomorphic Poisson
manifold.} A smooth real ( {\rm resp.} complex ) manifold $P$
 is called a {\it  Poisson manifold} if the algebra of smooth functions
 ${\cal A}=C^{\infty }(M,\R)$  ( {\rm resp.} $C^{\infty }(M,\C)$ )
 is equipped with a structure of Poisson
algebra over $\R $   ( {\rm resp.}  $\C$ ).

Let $P$ be a Poisson manifold of dimension $n.$ The Leibniz rule of
 $C^{\infty}(P)$ implies that the Poisson bracket $\{ , \}$ is
  a derivation in each argument, $\ie,$
$\{f,\cdot\}\in Der(C^{\infty}(M))$ is a vector field
on $P.$ By the skew-symmetry of $\{,\},$
 we can conclude that there exists a bi-vector $\pi_{P}
\in \Gamma(\W^{2} TP)$ such that $X_f := \{f,\cdot\}$ and
$\{f,g\}=\la\pi_{P}^{\sh}(df),dg\ra,$ here $\pi_{P}^{\sh}$ is the
bundle map defined by $\pi_{P}.$ The Jacobi identity just
corresponds to $[\pi_{P},\pi_{P}]=0,$
 here  $[,]$ is the {\it Schouten bracket} operation (\cf. [V, Chapter 1]).
 If $(P,\pi_P, J)$ is a holomorphic Poisson manifold with
  complex structure $J,$ clearly we have $\pi_{P}\in \Gamma(\W^{2} T^{1,0}P),$
  here  $T^{1,0}P$ is the holomorphic tangential space defined by $J.$
  A Poisson structure on a manifold $P$ defines a Lie algebroid structure on $T^*P.$ The Lie bracket on
  the space $\Gamma (T^*P)$ of  1-forms is given by
  $$\{\alpha,\beta\}= d\pi_P (\alpha,\beta)-\la\pi_P^{\sh}(\alpha),d\beta\ra + \la\pi_P^{\sh}(\beta),d\alpha\ra. \eqno{(2.1)}$$
  If $V\in \Gamma (TP)$ is a vector field, then
   $$\la V,\{\alpha,\beta\}\ra=(L_V\pi_P)(\alpha,\beta)-\la\pi_P^{\sh}(\alpha),d(i(V)\beta)\ra + \la\pi_P^{\sh}(\beta),d(i(V)\alpha)\ra. \eqno{(2.2)}$$

  If $\pi_P$ is a nondegenerate tensor,
the inverse bundle map
  $(\pi_P^{\sh})^{-1}$ thought of as a two form is a symplectic
 form, thus $P$ is the unique symplectic leaf. However in general case,
   $\pi_P$ may have varying ranks. Fix a point $p\in P,$ define
      $${\cal C}=\{V\in {T_p} P{\mid}  f\in C^{\infty}(P), {X_f} (p)=V\}.$$
      $\CN$ is a differential distribution on $P.$  By  Jacobi
      identity we have $[X_f, X_g]=X_{\{f, g\}}.$ So $\CN$
is an involutive distribution. In fact,
      $\CN$ is a completely integrable distribution (\cf.[V, Theorem 2.6]).
We denote
       the integral leaf through $p$ by $\fl_p.$ Since
 Hamiltonian
       flows defined by  Hamiltonian vector fields preserve the    Poisson
       brackets, the restriction of $\pi_P$ to $\fl_p$ has constant rank. In this way
        we know  that $\fl_p$ is a symplectic submanifold of $P$ and hence
        a symplectic leaf of $P.$  If $P$ is compact, let
        $\Ham (P)\subset \Diff (P)$ be the Hamiltonian diffeomorphism group with
        Lie algebra $\{X_f\mid f\in C^{\infty}(P)\}.$ Then the symplectic
        leaves are just the connected components of $\Ham (P)$-orbits in $P.$

        Let $r(p):=\dim  {\cal C} _p,$ then $r(p)$ is equal to the rank of $\pi_P (p).$
        Since $\pi_P$ is skew symmetric, $r$ is a even integer-valued function
        bounded by the dimension of $P.$  Choose a local coordinate chart
        $U\ni p=(x_{1},\cdots, x_{n})$ around $p.$ We can write under this
        coordinate that  $ \pi_P(p)=\pi_P^{ij}(p)\partial _{x_i }\W  \partial _{x_j}.$
        Fix a positive integer $k\leqslant n,$ for vectors
        ${\overrightarrow{\rho}}=(1\le \rho_1<\rho_2<\cdots <\rho_k \le n)$ and
         ${\overrightarrow{\delta}}=(1\le \delta_1 <\delta_2 <\cdots <\delta_k\le n),$ let
         $\pi^{{\overrightarrow{\rho}\overrightarrow{\delta}}}_P$
          be k-order
         sub-matrix of $\pi_P$ whose entries, lines indexed by
         $\overrightarrow{\rho}$ and  columns
         indexed by $\overrightarrow{\delta}.$ Define
         $$r_k (p)=\sum_{\mid{\overrightarrow{\rho}}\mid,\mid{\overrightarrow{\delta}}\mid=k}
         \mid{\rm det}\pi^{\overrightarrow{\rho}\overrightarrow{\delta}}_P(p)\mid ^2,
         \quad k=1,2,\cdots,n, $$
         here $|\overrightarrow{\rho}|=k$ means $\overrightarrow{\rho}$ is a vector with
         $n$-components.
         Let $P_l=\{x\in P\mid r(x)= l\}$ and
         $Q_l=\{x\in P\mid r(x)\le l\}.$ Denote `max' the maximum of $r(x).$
          Then

         $$P=\bigcup_{l=1}^{[n/2]} P_l=Q_{\max}.$$

        \bli{Proposition 1.1}  If $P$ is an holomorphic Poisson manifold,  then
         for  any positive integer less than $n/2,$
                 $P_{2l-1}$ is empty;
         if $Q_{2l}$ is non-empty,
         then the real codimension of the subset $Q_{2l-2}$ in
         $Q_{2l}$ is of at least 2.
         \eli
         \pf Clearly $P_{2l-1}$ is empty since the rank of a skew symmetric
         matrix is always even. We have

 $$
         Q_{2l}=\bigcup_{k=1}^{l} P_{2k}=
         \{x\in P\mid r_{2k+1}(x)=r_{2k+2}(x)=\cdots=r_{n}(x)=0\}\eqno{(2.3)}
$$
         and
$$
    Q_{2l-2}=\{x\in Q_{2l}\mid r_{2k}(x)=r_{2k-1}(x) =0\}\eqno{(2.4)}.
$$
  Note that all $Q_{2k}$ are analytic subset of $P.$
   Since  $Q_{2l}$ is non-empty,  $r_{2l}$ is a
    nonzero holomorphic function on $P.$ So the real
          codimension of  $Q_{2l-2}$ in $Q_{2l}$
 is of at least 2.  \qed

         \bli{Corollary 1.1} For holomorphic Poisson manifold $P,$
         if  $P_{2l}$ is non-empty,
         then  $P_{2l}$ is a relative open subset in
         $Q_{2l}.$   In particular,  $P_{\max}$ is an open connected density
       subset of $P$ if $P$ is a connected holomorphic Poisson manifold.
       \eli
             \pf  From equation (2.3) and (2.4), we know
         $Q_{2l}=P_{2l}\bigcup Q_{2l-2}.$  Note $Q_{2l}$ and $ Q_{2l-2}$ are closed,
    so  $P_{2l}$ is a relative open in  $Q_{2l}.$
Since $P_{\max}$   is the complement of the closed
          subset $Q_{\max-2}$ in $Q_{\max}$ and a real codimension two subset
  doesn't destroy the connectedness  of a manifold,
           $P_{\max}$ is an open connected density
         subset of $P.$     \qed

        \hskip 0.5cm  Any smooth Poisson manifold $P$ has the rank decomposition
         $P=\sum P_{2l},$ where $P_{2l}$ is foliated by symplectic leaves of dimension
         $2l.$ There is an open  subset $P_{\max}$ which is of maximal rank, but it
         does not necessary dense in $P.$  (Here is one such example.  The Poisson manifold
         $P=\R^2$ with Poisson  tensor $\pi_{P}= f(x,y)\pa_x \W\pa_y.$ Here $f$ is a smooth
         function such that
         $f(x,y)=0$ if $ y\geq 0, $ and $f(x,y)=x$ when $ y< 0.$ Note
         $P_{\max}=\{(x,y)\in R^{2}|x\not=0,y<0\}$ in this example). Since Hamiltonian
        flows preserve the rank of the Poisson
         bivector field, all $P_{2l}$ are
$\Ham (P)$-invariant subsets of $P$ (Assume $P$
is compact when $\Ham (P)$ is not appropriately
          defined). Denote the connected components
of $\Ham (P)$-orbits
         by ${\cal F}_i,$ where $i$ is indexed by
 an index set denoted by $\L_l.$ Then we have a symplectic
 decomposition
         $$P=\sum_{\lambda\in \L _l,  1\leq l\leq [n/2]}
{\cal F}^{\lambda}_l.$$
         The symplectic leaves in $P_{\max}$
are all open. So the $\Ham (P)$-action has at least one open orbit.
          If $(P, \pi_P, J)$ is a holomorphic symplectic
 manifold, the Hamiltonian vector field of any holomorphic
 function is holomorphic. So, all
          symplectic leaves are complex
symplectic submanifolds and the restriction of $\pi_P$ to the
symplectic leaf ${\cal F}_i $ defines a symplectic
          form ${\omega}^{\lambda}_l$
          satisfying ${\omega}^{\lambda}_l(J\cdot,J\cdot)=
          {\omega}^{\lambda}_l(\cdot,\cdot).$
          Hence  ${\cal F}_i $ are
          pseudo-K\"ahler submanifolds. Thus we have

         \bli{Proposition 1.2} Let $P$ be a  holomorphic Poisson
         manifold.  Then we have a symplectic decomposition
         $$P=\sum_{\lambda\in \L _l,  1\leq l\leq [n/2]}
{\cal F}^{\lambda}_l.$$
          Where each symplectic leaf ${\cal F}^{\lambda}_l$ is
         a pseudo-K\"ahler submanifold. Moreover $P$ has at
        least one open  symplectic leaf.
       \eli

\sect{\S3. \ Poisson Lie group}
  Before giving the definition of a Poisson Lie group,
   we need some notations.
   Let ${\cal A, B} $ be Poisson algebras over $\F.$
A homomorphism $f: {\cal A}\lw {\cal B}$ is
 called a {\it Poisson homomorphism} if it is bracket preserved, $\ie,$
 $f(\{a,b\})=\{f(a),f(b)\}.$  A smooth
 map between smooth Poisson manifolds  $M$ and $N$ is called a
{\it Poisson map} if the pull back map $f^{*}: C^{\infty}(N)\lw
C^{\infty}(M) $ is a Poisson homomorphism.  Given two Poisson
manifolds $(M, \pi_{M})$ and $(N,\pi_{N}),$ there is a
 unique Poisson structure, denoted by
$\pi_{M}\oplus \pi_{N},$
 on the product manifold $M\times N,$ such that
the projections from $M\times N$ to each factor $M$ and $N$ are
both Poisson maps. We called it the {\it product Poisson
structure.} For algebraic varieties (analytic spaces), we have
similar definitions of  Poisson morphisms
 and Poisson structure on Product varieties (spaces).

A Lie (resp. algebraic) group which is at the
same time a Poisson manifold is called a
{\it Poisson Lie (resp. algebraic) group} if
 the multiplicative map
     $\mu :G\times G \lw G,  (h_{1},h_{2}) \lo  h_{1}h_{2}$
  is a Poisson map, where $G\times G$
is equipped with the product Poisson structure.

       Let $r_{g}$ and $l_{g}$ respectively denote
 the left and right translation in $G$ by $g.$
Then the Poisson bivector field  $\pi_{G}$ of a Poisson Lie group satisfies
$$\pi_{G}(gh)=l_{g}\pi_{G}(h)+r_{h}\pi_{G}(g),        \forall   g, h\in G,$$
and is called a {\it multiplicative} bivector field. If both $g,h$
are specialized to be the unit $e,$
 then $\pi_{G}(e)=0.$ So a Poisson Lie groups
can never be a symplectic manifold.
 Let $\pi^r  (g)=r_{g^{-1}} \pi_G (g).$ Then
 $$
 \pi^r (gh)=\pi^r (g)+Ad_{g}\pi^r (h)\eqno{(3.1)}
 $$
    The linearization of $\pi^r$ at the unit  is a $\W^2 \fg$-valued
 $1$-cocycle on $\fg$ relative to the adjoit
action of $\fg$ on $\wedge^{2}\fg,$  that is $ \delta=d_{e}\pi^r
\in H^{1}(\fg, \wedge^{2}\fg).$
  If  $\delta $  is a 1-coboundary,  then $G$
  is called a  {\it{ coboundary  Poisson-Lie   group}}.
 The dual map  $\delta ^{*}$
 defines a  Lie algebra structure
   on the dual space  $\fg^{*}$  of  $\fg.$
There is a Lie algebra structure on $\fg\oplus\fg$ such that both
$\fg$ and $\fh^*$ are Lie subalgebras.
    $(\fg\oplus\fg,\fg,\fg^{*})$ is called the
{\it Lie bialgebra} of $(G,\pi_G).$
    There is a 1-1 correspondence
   between connected and simplely connected
 Poisson Lie groups and Lie bialgebras (\cf. [LW]).
The connected and simply connected group $G^{*}$
   with   Lie algebra $\fg^*$ is also a Poisson Lie
 group, called the {\it dual} Poisson Lie group of $G.$
      A Lie group $D(G)$ with Lie algebra
$\fg\oplus \fg^*$ is called a {\it double} for
 the Poisson Lie group $G$ if
   the multiplication map
$G^{*}\times G\lw D(G),(g,h)\lo gh$ is a diffeomorphism.
In this case, the left-action of $D(G)$ on itself induces an
   action  on $G^{*}.$ Its restriction
$G\times G^* \lw G^*, (g,h)\lo D^l (g)$ i
s called the {\it left
dressing action} of $G$  on $G^{*}.$
   If $G$ has trivial Poisson Lie structure,
$G^*=\fg^*$ has Lie Poisson structure,
 and the left drssing action is
 just coadjoint action of $G$
   on $G^*.$ The definition of
 right dressing action is similar (\cf. [Lu]).

    If $G$  is semisimple,  the  Whitehead Theorem says
 that   both $H^{1}(\fg,V)$ and $H^{2}(\fg,V)$ are vanishing for any
$\fg$-module $V.$  Thus every connected
semi-simple Poisson Lie group is coboundary
Poisson Lie group. The $1$-cocycle is given by
  $\delta(X)=ad_{X}\L$ for some $\L\in\wedge^{2}\fg.$
 By integration  we get a 2-cycle  $\pi_G (g) =l_{g}\L-r_{g}\L$
on the group $G.$ Clearly $\pi_G$ is a Poisson bivector  if and
only if $[\L,\L]\in (\wedge\fg)^{\rm inv}.$   Such $\L$ is called
a {\bf r}-{\it matrix}.

       If    $G$ is an abelian group, then $g\lw \pi^r (g)$  is a Lie group
homomorphism from $G$ to abelian group $\W^2 \fg.$  So on the compact
       abelian group there is no nontrivial  Poisson Lie group structure.

      If  $G=T^m \times \R^n,$ a direct product of torus group and vector group.
         Let
         $$\overline{\pi}^r (\theta_{1},\cdots,\theta_{m},x_{1},\cdots,x_{n})
={\pi}^r (e^{i\theta_{1}},\cdots,e^{i\theta_{m}},x_{1},\cdots,x_{n})$$
         be the lift of $\pi^r$ on the universal covering space $\R^{m+n}.$
 By the multiplicative condition (3.1), we know $\overline{\pi}^r_{ij}$ is a   linear
         function with respected to the variables $\theta_{\rho}$ and $x_{\delta}$
         and periodic with respected to $\theta_{\rho}.$
         So the Poisson bivector  has the form
$$\overline {\pi}^r(u_{1},\cdots,u_{m+n})=
\sum^{m+n}_{i,j,k=1}C^k_{ij}u_{k}{\partial}_{u_{i}}\W\partial_{u_{j}},$$
 where
the $C^k_{ij}$ are the structure constants of a $m+n$-dimensional
Lie algebra
 and $C^{k}_{ij}=0$ for $k\le m;$ $u_{i}=\theta_{i}$
for $1\le i\le m$ and $u_{m+i}=x_{i}$ for $1\le i\le n.$
  In particular, if $G$ is the vector group
$\R^m,$ the Poisson Lie structure on $G$
is linear and called a {\it Lie Poisson} structure.
  For multiplicative abelian group  $G=(\R_{+})^n$
  or  $(\C^{\star})^{n}.$ the solution of (3.1) is of form
   $$
 \pi^r_{\mu\nu} (z_{1},\cdots, z_{n})=\sum^n_{\delta=1} C^\delta_{\mu\nu}{\rm ln}z_{\delta}, \eqno{(3.2)}$$
 $\pi_G$ is a Poisson structure iff $C^\delta_{\mu\nu}= -
C^\delta_{\nu\mu}\in\F  $ for $ 1\le \mu,\nu,\delta\le n$ and
 $$
 \sum_{\mu,\nu=1}^n[  (C^{\mu} _{\rho \delta} C^{\nu}
_{\delta\gamma}+C^{\mu} _{\rho \gamma} C^{\nu}
_{\delta\gamma}){\rm ln}z_{\mu}{\rm ln}z_{\nu}+
 C^{\mu} _{\rho \nu} C^{\nu} _{ \delta\gamma}
{\rm ln}z_{\mu}] +{\rm c.p.}(\rho,\delta,\gamma)=0
 \eqno{(3.3)}
 $$
for $1\le \rho,\delta,\gamma\le n,$   here ${\rm c.p.}
(\rho,\delta,\gamma)$ means cyclic permutation with respected to
$\rho,\delta$ and $\gamma.$
 After some simple calculations, we know (3.3)  is equivalent to
 $$
 \sum _{\nu=1}^{n}
 C^{\mu} _{\rho \nu} C^{\nu} _{\delta\gamma}
 +C^{\mu} _{\delta \nu} C^{\nu} _{Z\rho}
+C^{\mu} _{\gamma \nu} C^{\nu} _{\rho\delta}=0,\quad
 1\le \mu,\rho,\delta,\gamma\le n. \eqno{(3.4)}
  $$
  So $C^{\mu}_{\rho\delta}$ are also the
structure constants of a $n$-dimensional Lie algebra.
  The Poisson Lie group structure  on $G$ is given by
    $$\pi_{G}(z_{1},\cdots, z_{n})=\sum^n_{ \delta,\mu, \nu=1}
C^\delta_{\mu\nu}z_{\mu}z_{\nu}{\rm ln}z_{\delta}
\partial_{\mu}\wedge \partial _{\nu}$$

    \bli{Example 3.3}  $\pi=x(ae^{i(\theta_{1}+\theta_{2})}\partial_{\theta_{1}}\wedge
\partial_{\theta_{2}}+be^{i\theta_{1}}\partial_{\theta_{1}}\wedge\partial_x
    +ce^{i\theta_{2}}\partial_{\theta_{2}}\wedge\partial_x)$
    is a Poisson Lie group structure on $T^2 \times \R $  for any
    $a,b,c, \theta_{1},\theta_{2} \in \R.$
    \eli

    Belavin and Drinfeld classified all Poisson Lie structures
on complex simple Lie group [BD], the
    classification of Poisson-Lie structure  on  compact
connected Lie group was carried out by Levendorskii and Soibelman
in [LS].
   The classification of Poisson Lie structure on
 non-semisimple non-abelian group is more involved.
   In [M], [Z],  [BKM], where the cases
     2-dimention affine group, Poincar\'e group,
Galilei group, respectively,  were solved completely.

\sect{\S4 Poisson actions and momentum mappings} \hskip\parindent
An action of a Poisson Lie group $(G,\pi_{G})$ on a Poisson
manifold $(P,\pi_{P})$ is called a {\it Poisson action} if the
action map $\sigma : G\times P\lw P$ is a Poisson map, where
$G\times P$ is equipped with product Poisson structure. Let $g\in
G$ and $x\in P,$ denote by $\si_{g}:P\lw P$ and by $\si_{x}:G\lw
P$ the maps
$$\si_{g}:x\lo \si(g,x)=gx,\quad   \si_{x}:g\lo \si(g,x)=gx. $$
  $\si$ is a {\it Poisson action} iff
$$\pi_{P}(gx)=\si_{g_{*}}\pi_{P}(x)+\si_{x*}\pi_{G}(g).$$
Denote the infinitesimal action by $\l:\fg\lw \G(TP),$
equivalently, $\l(X)(x)=\si_{x*}X.$
 For any $f\in C^{\infty}(P),$ let $\xi_{f}\in \fg^{*}$
  be the covector defined by  $$\langle\xi_{f}(p),X\rangle=\frac{d}{dt}
\Big|_{t=0}f(e^{tX}p)=\langle df(p), \l(X)\rangle.\eqno{(4.1)}$$
  If $G$ is connected, $\sigma$ is a Poisson action iff
  $$\l(X)(\{f,g\})= \{\l(X)(f),g\}+
\{f, \l(X)(g)\}+\langle
X,[\xi_{f},\xi_{g}]_{*}\rangle,\eqno{(4.2)}$$
Clearly, if $G$ has trivial Poisson structure,
 in this case $G^{*}=\fg^*$ is an abelian Poisson Lie group, the
$G$-action  preserves the Poisson structure of $P.$ Conversely, we have

\bli{Proposition 4.1 } Let $\sigma: G \times P \lw P$ be  a
Poisson action. Then it preserves the Poisson structure of $P$ if
and only if for every $p\in P,$ the annihilator $\fg^{\circ}_{p}$
 of the isotropy subalgebra $\fg_{p}$ at $p$ is abelian. In particular,
 if the action is locally free, then $\sigma$ is a Poisson structure
 preserved action if and only if $G$ is a trivial Poisson Lie group.
 \eli

\pf     $\sigma$
  preserves the Poisson structure of $P$ if and only if
$$\l(X)(\{f,g\})=
\{\l(X)(f),g\}+\{f, \l(X)(g)\},\eqno{(4.3)}$$
Define the map $A:d(C^{\infty}(P))\lw \fg^{*},df\lo \xi_{f}.$
Then by (4.1), the image of $A$ is the annihilator of
$\fg_{p}=\{X|\l(X)(p)=0\}.$
By (4.2) and (4.3), we know  the action preserves
 the Poisson structure of $P$ if
  and only if for every $p\in P,$ the annihilator
$\fg^{\circ}_{p}$ of the isotropy subalgebra $\fg_{p}$
 is abelian. So the linearization of Poisson structure
of $G$ is zero, which means the Poisson Lie structure
 of $G$ is trivial.
\qed

    Poisson action $\si:G\times P\lw P$ is said to be
    {\it tangential } if every infinitesimal vector
    field $\l(X)$ is tangent to the symplectic leaf of $P.$
       A smooth map $m:P\lw G^*,$  is called a {\it momentum map}, if $\l(X)=\pi^{\sh}(J{^*} X^l),$ here
       $X^l$ is the left invariant 1-form on
       $G^*$ whose value at the identity is
        $X\in \fg=(\fg^*)^*.$ A momentum map is called equivariant if
         it intertwines
       with the action of $\si$ and the left dressing
        action of $G$ on $G^*.$  A Poisson action with
         momentum map is clearly tangential action.
        Dressing actions are tangential Poisson actions
        (\cf. [STS, LW]). The infinitesimal vector field
          of left  dressing action of $G$ on $G^{*}$ is given by
          $$\d^l(X)=\pi_{G^*}^{\sh}(X^l). $$
         Clearly by the definition of dressing action, the dressing orbits sweep out all
         symplectic leaves of $G^{*}$ and the identity map
          is an equivariant  moment map.

In general case if there exists a momentum map, it is not
necessary equivariant.  Let $m$ be a momentum map  which is not
necessary eqivariant, define a map
$${\Sigma}:G \times P \lw  G^{*},(g,x) \lo (D^{l}(g)m(x))^{-1}\cdot m(gx).$$
 $\Sigma$ measures the equivariant properties  of m,
 it is  equivariant if and only if the image of
  $\Sigma$ is a single point.
We want to give an infinitesimal description of $\Sigma.$ For that
we fix $x \in P ,$  denote ${\Sigma}^{x}: G \lw G^{*},
 g \lo {\Sigma}(g,x)$ and
  ${\Gamma}_{{X},{Y}}(x)={\langle}d_{e}{\Sigma}^{x}
  (X),Y{\rangle},$  where $X,Y\in \fg.$  Then we have

  \bli{Proposition 4.2}
{\rm (i)}$$\Gamma_{X,Y}=m^* (\pi_{G^*}(X^{l},Y^{l})) -\pi_{P}(m^*
X^l , m^* Y^l); \eqno{(4.4)}.$$ {\rm
(ii)}$$\pi_{P}^{\sh}(d{\Gamma}_{X,Y})= \pi_{P}^{\sh}( m^*( i(\d^l
(X)) dY^{l})-i(\d^l (Y)) dX^l )- i(\l(X))dm^*Y^l+i(\l (Y))d m^*
X^{l})\\\eqno{(4.5)} .$$

{\rm (iii)}  For $Z \in \fg,$ the differential of
${\Gamma}_{{X},{Y}}$  in the direction of $\l(Z)$ is
$$\begin{array}{rcl}
\l(Z)({\Gamma}_{X,Y})&=&{\langle}\l(Z), m^*
(\{X^{l},Y^{l}\}_{\pi_{G^{*}}})
-\{ m^* X^l, m^*Y^l\}_{\pi_{P}}{\rangle}\\
&&  +
  dX^{l}( m_{*}(\l(Y))-\d^l(Y),m_{*}(\l(Z)))\\
  && + dY^{l}( \d^l(X)- m_{*}(\l(X)),m_{*}(\l(Z)))
 \end{array}
\eqno{(4.6)}$$

{\rm (iv)}
 Define the  map
 ${\Gamma}: \fg\times \fg \lw C^{\infty}(P),(X,Y) \lo {\Gamma}_{X,Y}.$
 Then ${\Gamma}$ is antisymmetric and bilinear and
 $$\begin{array}{rcl}
 &&\hskip 1.0cm (d_{*} {\Gamma}+\frac{1}{2}  [{\Gamma},{\Gamma}])(X,Y,Z)\\
 &&=\{{\Gamma}([X,Y],Z)
 -{\langle}Z,[(m_{*}(\l(X)))^{l}-(\d^l (X))^{l},(m_{*}(\l(Y)))^{l}-(\d^l (Y))^{l}]_{*}{\rangle}\}\\
&& \hskip 0.5cm +{\rm c.p.}(X,Y,Z)\\
 &&=\{{\Gamma}([X,Y],Z)
 -[(L_{\l(Z)}\pi_{P})(m^{*}X^{l},m^{*}Y^{l})+
  (L_{\d^l (Z)}\pi_{G^{*}})(X^{l},Y^{l})]\\
&&\hskip 0.5 cm  +[ (L_{\overline{Z}}\pi_{G})(e)
((\d^l (X))^{l},(m_{*}(\l(Y)))^{l})+
 (L_{\overline{Z}}\pi_{G})(e)
((m_{*}(\l(X)))^{l}, (\d^l (Y))^{l})]\}\\
&&\hskip 0.5cm +{\rm c.p.}(X,Y,Z);
 \end{array}
\eqno{(4.7)}$$ here  $d:\Gamma(\bigwedge^* T^*
P)\lw\Gamma(\bigwedge^{*+1} T^* P)$ is usual de Rahm differential
on forms decided by the Lie brackets of vector fields on $P;$ and
$ d_{*}: {\rm Hom}(\bigwedge^* \fg, C^{\infty}(P))\lw {\rm
Hom}(\bigwedge^{*+1} \fg, C^{\infty}(P))$
 is the de Rham differential decided by the Lie bracket of $\fg
 ^{*}$  which we denoted by  $[,]_{*};$
  and ${\overline {Z}}$  is any vector field on $G$ whose value at the unit of $G$ is $Z.$
  \eli
  \pf (i)
${{\Sigma}}^{x}$ maps the unit of $G$ to that of $G^{*},$ so there
is  an induced map $d_{e}{{\Sigma} ^{x}}: \fg \lw {\fg}^{*}.$ And,
 $$\begin{array}{rcl}
 {\Gamma}_{X,Y}&=&{\langle}d_{e}{\Sigma}^{x}({X}),{Y}{\rangle}\\
 &=&{\langle}(l_{{m(x)}^{-1}})_{*}(m_{*}(\l(X))(x))-
 (l_{{m(x)}^{-1}})_{*}(\d^l (X)(m(x)),Y{\rangle}\\
 &=&{\langle}m_{*} {X}_{P}(x)-\d^l (X)(m(x)),{Y}^{l}{\rangle}\\
 &=&{\langle}m_{*} {X}_{P},{Y}^{l}{\rangle}+{\langle}\pi_{G^{*}}^{\sh}  ({{X}^{l}})\circ m,{Y}^{l}{\rangle}\\
 &=&{\langle}\pi_{G^{*}}^{\sh}({X}^{l}),{Y}^{l}{\rangle}\circ m-
{\langle}m_{*}(\pi_{P}^{\sh}(m^{*}{X}^{l})),{Y}^{l}{\rangle} \\
 &=&m^{*}\pi_{G^{*}}({X}^{l},{Y}^{l})-
 \pi_{P}(m^{*}{X^l},m^{*}{Y^l}).
 \end{array}$$

(ii)  Using the defintion (2.1), we have,
       $$\begin{array}{rcl}
d{\Gamma}_{{X},{Y}}
&=& m^{*}(\{{X}^{l},{Y}^{l}\}_{\pi_{G^{*}}}-
i({\pi_{G^{*}}^{\sh}({X}^{l})}) {Y}^{l}+
i({\pi_{G^{*}}^{\sh}({Y}^{l})}) {X}^{l})\\
&&\hskip 0.07cm
-\{m^{*}{X^l},m^{*}{Y^l}\}_{\pi_{P}}+
i({\pi_{P}^{\sh}(m^{*}{X^l})}) dm^{*} ({Y}^{l})-
i({\pi_{P}^{\sh}(m^{*}{Y^l})}) dm^{*} ({X}^{l})\\
&=&m^{*}(\{{X}^{l},{Y}^{l}\}_{\pi_{G^{*}}})
-\{m^{*}{X^l},m^{*}{Y^l}\}_{\pi_{P}}\\
&&\hskip 0.07cm -i{(\l(X))} d m^{*} ({Y}^{l})+
i{(\l(Y))} dm^{*} ({X}^{l})
- m^{*}(-i({\d^l (X)}) d{Y}^{l}+
i({\d^l (Y)}) d{X}^{l})\\
  \end{array}.$$
  Using the following identity
   $$\begin{array}{rcl}
&&\pi_{P}^{\sh}(  m^{*}(\{{X}^{l},{Y}^{l}\}_{\pi_{G^{*}}})
-\pi_{P}^{\sh}\{m^{*}{X^l},m^{*}{Y^l}\}_{\pi_{P}})\\
&=&\pi_{P}^{\sh}( m^{*}([X,Y]^{l})-
[\pi_{P}^{\sh}(m^{*}{X^l}),\pi_{P}^{\sh}(m^{*}{Y^l})]\\
&=&\l([{X},{Y}])+[\l(X),\l(Y)] =0
     \end{array}$$
to cancel the first two terms in the expression of $d\Gamma_{X,
Y},$ then applying map $\pi^{\sh}$ to it  we get (4.5).

(iii)  Firt note that $\l(Z)({\Gamma}_{X,Y})=\langle \l(Z),
d{\Gamma}_{X, Y} \rangle,$ since $$
\begin{array}{rcl}
&&i(\l(Z))(m^{*}(i(\d^l (X)) d{Y}^{l})
- i(\l(X)) dm^{*} ({Y}^{l})+
i(\l(Y)) dm^{*} ({X}^{l})
-m^{*}(i(\d^l (Y)) d{X}^{l}))\\
=&& dX^{l}(m_{*}(\l(Y))-\d^l (Y),m_{*}(\l(Z)))+
 dY^{l}(\d^l (X)-m_{*}(\l(X)),m_{*}(\l(Z))),
 \end{array}
$$
and the calculations in (ii) we can easily get (4.6).

 (iv)Using
$$\frac{1}{2}[{\Gamma},{\Gamma}] (X,Y,Z)=
{\langle}Z,{\Gamma}[X,Y]_{{\Gamma}}-[{\Gamma} X,{\Gamma} Y]_{*}{\rangle},
$$and
$${\langle}Z,{\Gamma}[X,Y]_{{\Gamma}}{\rangle}
=-{\langle}{\Gamma}Z,ad_{{\Gamma} X}^{*}Y-ad_{{\Gamma}Y}^{*}X{\rangle}
=-{\langle}Y,[{\Gamma}Z,{\Gamma}X]_{*}{\rangle}-{\langle}X,[{\Gamma}Y,{\Gamma}Z]_{*}{\rangle},$$
we know
$$  \begin{array}{rcl}
&&\hskip 0.07cm (d_{*} {\Gamma}+\frac{1}{2} [{\Gamma},{\Gamma}])(X,Y,Z)\\
 && ={\langle}Z,{\Gamma}[X,Y]-[{\Gamma} X,{\Gamma} Y]_{*}{\rangle} + {\rm c.p.}(X,Y,Z)\\
 &&={\langle}Z,(m_{*}[X,Y]_{P})^{l}-([X,Y]_{G^{*}})^{l}{\rangle}
 -[(m_{*}(\l(X)))^{l}-(\d^l (X))^{l},(m_{*}(\l(Y)))^{l}-(Y_{G_{*}})^{l}]_{*}{\rangle}\\
 &&\hskip 0.5cm +{\rm c.p.}(X,Y,Z)\\
 &&=\{{\Gamma}_{[X,Y],Z}
 -{\langle}Z,[(m_{*}(\l(X)))^{l}-(\d^l (X))^{l},(m_{*}(\l(Y)))^{l}-(Y_{G_{*}})^{l}]_{*}{\rangle}\}\\
 &&\hskip 0.5cm +{\rm c.p.}(X,Y,Z).
  \end{array}$$
 Since $G \times P \lw P$  and  $G \times G^{*} \lw  G^{*}$ are Poisson actions,
by Theorem 4.7 of {[Lu]},
$$ (L_{\l(Z)} \pi_{P})(m^{*}X^{l},m^{*}Y^{l})=
{\langle}Z,[(m_{*}(\l(X)))^{l},(m_{*}(\l(Y)))^{l}]_{*}{\rangle};
\eqno{(4.8)}
$$
$$(L_{\d^l(Z)} \pi_{G^{*}}) (X^{l},Y^{l})
 ={\langle}Z,[(\d^l (X))^{l},
(\d^l (Y))^{l}]_{*}{\rangle};
\eqno{(4.9)}$$
 And by the definition of $[,]_{*},$  we have,
$$(L_{\overline{Z}}\pi_{G})(e) ((\d^l (X))^{l},(m_{*}(\l(Y)))^{l})
={\langle}Z,[(\d^l (X))^{l},(m_{*}(\l(Y)))^{l}]_{*}{\rangle};\eqno{(4.10)}$$
$$(L_{\overline{Z}}\pi_{G})(e)(({\rm
 m}_{*}(\l(X)))^{l}, (\d^l (Y))^{l})={\langle}Z, [{\rm
 m}_{*}(\l(X)))^{l}, (\d^l (Y))^{l}]_{*}{\rangle}.
\eqno{(4.11)}$$
By (4.8-4.11), we have the last identity of (4.7).
\qed

\vskip0.5cm
     If $G$ is connected,
    the momentum  map $ m :P \lw G^{*}$ is equivariant  if and only
    if $m$ is a Poisson map by (ii) of Proposition 4.2.
 If $P$ is a symplectic manifold, and $x_{0}\in P$ such that
 $m(x_{0})=e^{*}$ is the unit of $G^{*}.$ Let $\L=m_{*}\pi_{x_{0}}\in
  \bigwedge ^{2}\fg^{*},$
 and $\pi_{\tw}=\pi_{G^{*}}+\L^{r}.$  Then $\pi_{\tw}$ is still a Poisson
 structure of $G^*.$
 It is proved in [Lu] that $m:(P,\pi_{P})\lw (G^{*},\pi_{\tw})$ is
  an equivariant momentum map. Which means
  $m^* ((\pi_{G^*}+\L^r)(X^{l},Y^{l})) -\pi_{P}(m^* X^l , m^* Y^l)=0.$ So
  $$\Gamma(X,Y)(x)=(Ad(m(x)^{-1})\L)(X,Y).$$
  Since $\pi_{\tw}$ is a Poisson structure, we have $d_{*}\L+\frac{1}{2}[\L,\L]=0.$
  But we even don't know when $d_{*}\G+\frac{1}{2}[\G,\G]=0$ on
  the symplectic leaves of $P.$

From now on we suppose $G$ is a trivial Poisson Lie group. Then
$\pi_{G}=0$ and $G^{*}=\fg^{*}.$
 Denote the set of Carsimir functions on $P$ by
 ${\rm Car}(P).$ In this case ${\Gamma}_{X,Y}\in {\rm Car} (P)$ is a Carsimir
 function and
  $$\Gamma_{X, Y}(x)=  m([{X},{Y}])(x)-
 \{m({X}),m({Y})\}_{\pi_{P}}(x),$$
where we denote the $X$-component of $m(x)$ by $m(X)\in
C^{\infty}(P).$
 By definition the Poisson structure of $P$ is
preserved under the action of $G$ we have $L_{\l(Z)}\pi_{P}=
L_{\d^l (Z)}\pi_{G^{*}}
 =L_{\overline{Z}}\pi_{G}=0.$  So
$\frac{1}{2} [{\Gamma},{\Gamma}](X,Y,Z)=0$ and
     (4.7) is reduced to
   $$\begin{array}{rcl}
&&\hskip 0.07 cm d_{*} {\Gamma}(X,Y,Z)\\
& =&({\{ }m([{X},{Y}],m({Z}){\}}_{\pi_{P}}- m([[{X},{Y}],{Z}])) +{\rm c.p.}({X},{Y},{Z}) \\
& =&({\{} \{ m({X}),m({Y}) {\}}_{\pi_{P}} -{\Gamma}_{{X},{Y}},
m({Z}){\}}_{\pi_{P}}- m([[{X},{Y}],{Z}]) )+{\rm c.p.}({X},{Y},{Z})  \\
& =&({\{} \{ m({X}),m({Y}) {\}}_{\pi_{P}} ,m({Z}){\}}_{\pi_{P}}- m([[{X},{Y}],{Z}]) )
+{\rm c.p.}({X},{Y},{Z}) \\
& =& 0.\\
\end{array}$$
 Hence ${\Gamma}:\fg \times \fg \lw {\rm Car}(P) $ is a ${\rm Car}
 (P)$-valued
 2-cocycle of  $\fg$ relative  to the trivial representation of
  $\fg$ on ${\rm Car}(P),$ equivalently
 $[{\Gamma}] \in H^{2}(\fg, {\rm Car} (P) ) \cong  H^{2}(\fg) \otimes {\rm Car}(P).$   In the case that
 $H^{1}(\fg, {\rm Car} (P))=H^{2}(\fg, {\rm Car} (P))=0,$ in particular, if $\fg$
is semisimple, this action has  a coadjoit equivariant momentum map.

 Let ${\Gamma}_{g,{Y}}(x)={\langle}{\Sigma}(g,x),{Y}{\rangle}.$
 If $G$ is connected then ${\Gamma}_{g,{Y}}$
   is also a Carsimir  function on $P.$ Let  $ L(\fg,{\rm Car}
    (P)) \cong \fg^{*} \otimes {\rm Car}(P)$ denote
   the linear maps from  $\fg$ to   ${\rm  Car} (P).$
   Define a map  $\Psi $ from $G$ to $L(\fg,{\rm Car} (P))$ by
   $\Psi(g)({Y})={\Gamma}_{g,{Y}}.$  Then
 $$ (\Psi (gh)({Y}))(x)
 =(\Psi(g)({Y}))(hx)+(\Psi(h)(Ad_{g^{-1}}{Y}))(x)$$
 Because the action is tangential, $hx$  and $x$
  lie in the same symplectic leaf
 and  ${\Gamma}_{g,{Y}}$ is a Carsimir function,
  we have $ (\Psi(g)({Y}))(hx)=(\Psi(g)({Y}))(x).$ Hence
 $\Psi (gh)=\Psi (g)+Ad^{*}_{g^{-1}} \Psi (h), $ and $ \Psi$ is a
  $ L(\fg,{\rm  Car} (P))$-value 1-cocycle  on
 group $G,$ it  represents a  cohomology class
 in  $ H^{1}(G, L(\fg,{\rm Car} (P)).$
So  $G$ admits a coadjoint equivariant momentum map
  if and only if $\Psi$ is cohomologeous to zero. In particular,
    if $G$ is compact, then it admits an equivariant momentum map.

       \sect{\S5. \  Almost homogeneous tangential Poisson actions}
Let $G$ be a Lie group and $P$ a $G$-manifold. $P$ is called
$G$-{\it almost homogeneous} if $G$ has only one open orbit in
$P.$  Note here our definition of almost homogeneous is a little
different to Huckleberry's in [HO] where it is defined for complex
space. However ours definition can apply  to real case. If $G$ is
a Poisson Li group and $P$ is a $G$-almost homogeneous Poisson $G$
space, then $P$ is an {\it almost homogeneous Poisson
$G$-manifold.}
   If $P$ is an almost  homogeneous Poisson
manifold under the tangential action of $G,$
    then each  open $G$-orbit must be a symplectic
leaf of $P,$ since the $\fg$-vector fields, which at the same
       time are Hamiltonian, generate the whole tangential space.
    So, if an almost homogeneous Poisson $G$-action has an equivariant momentum map,
     then all open $G$-orbits in $P$ are symplectic covering spaces of left
    dressing orbits
     in $G^{*}.$  In particular, if $P$ is  homogeneous,
$P={\cal O}_{p}=G\cdot p=G/G_p$ for some $p\in P$
    is a symplectic covering space of the $G$-dressing
orbit ${\cal O}_{u}=G/G_u\subset G^*,$ here $ u=m(p)$ and
 $m:{\cal O}_{p}\lw {\cal O}_{u}$
   is a homogeneous fibration whose fibers  are isomorphic to $G_u/G_p.$

\bli{Proposition 5.1} Let $\sigma: G\times P\lw P$ be a Poisson
action on a symplectic manifold $P$ and ${\rm m}:P\lw G^{*}$ an
equivariant momentum map. Then
$${\rm Ker ({m_{p}}_{*}})=(T_{p} {\cal O}_{p})^{\omega};
\hskip 0.7 cm {\rm Im ({m_{p}}_{*})}=[ (\fg _{p}^l)]^{\circ},
\eqno{(5.1)}$$ here the upscript $\omega$ denotes the
 annihilator
operation with respect to symplectic structure of $P.$ \eli \pf
$\forall f \in C^{\infty}(P),$
$${\langle}(\l(X))(p),df{\rangle}=
{\langle}\pi_{P}^{\sh}(m^{*}X^{l}), df{\rangle}
=-{\langle}m^{*}X^{l},\pi_{P}^{\sh}(df){\rangle}=
-{\langle}X^{l}, m_{x*}(\pi_{P}^{\sh}(df){\rangle},
\eqno{(5.2)}$$ so
$$\omega((\l(X))(p),\pi_{P}^{\sh}(df))
 ={\langle}X^{l}, m_{*}(\pi_{P}^{\sh}(df)){\rangle},
\eqno{(5.3)}$$
Let $v=\pi_{P}^{\sh}(df) \in T_{p} P,$  then $v$ runs through $T_{p}P$
since  $P$ is a symplectic manifold,
$$\omega((\l(X))(p),v) ={\langle}X^{l},
m_{*}v{\rangle}, \forall v \in T_{p}P.$$
So ${\rm Ker ({m_{p}}_{*}})=\{(\l(X))(p)|X \in \fg\}^{\omega}=
(T_{x} {\cal O}_{p})^{\omega},$ and
 $ {(\rm Im({m_{p}}_{*}))}^{\circ}=
 \{X^{l}|(\l(X))(p)=0\}=\{X^{l}|X\in
\fg _{p}\},$ the later means ${\rm Im ({m_{p}}_{*})}=
[(L_{u}^{-1})^{*} (\fg _{p})]^{\circ}.$\qed

Now assume $G$ is a compact Poisson Lie group. Let
$$P_{\max}=\{x\in m(P)\mid \dim G_{m(x)} \le \dim G_{m(y)}, \forall
y\in P \}.$$

\bli{Lemma 5.1} The set $P_{\max}$  is open in $P.$  \eli
 \pf
Since $G$ is a compact Lie group, the action $G\times G^{*}\lw
G^{*}$ is a proper action. Thus we can apply the slice  theorem,
it follows that for any $u\in G^*$ there is a $G$-invariant
neighborhood $U\ni u$ with $\dim G_{u} \le \dim G_{v}$ for all
$v\in U.$ The result follows from the equivariance and continuity
of $m.$ \qed

       \bli{Proposition 5.2}
        For any $p\in P_{\max},$   we have $[G_u,G_u]\subset G_p,$ where $u=m(p).$
       \eli

\pf  Let $\gamma:[0,\epsilon]\lw P_{\max}$ be a smooth curve with
$\gamma(0)=p$ and $\gamma'(0)=V\in T_p P.$
  Note that $\dim  \fg_{\scriptscriptstyle   m(\gamma(t))}$
is a locally constant  near $\gamma(t),$ hence for any $X, Y\in
\fg_u,$ we can find smooth curves $X_t,$  $ Y_t
\in\fg_{\scriptscriptstyle   m(\gamma(t))}.$ So
$$\pi_{G^*}(m(\gamma(t)))(X^l_t,Y^l_t)= \la\pi^{\sh}_{G^*}(X^l_t),Y^l_t
\ra(m(\gamma(t))=\la d^{l}(X)(m(\gamma(t)),Y^l_t\ra =0,$$
Thus $$\pi^l_{G^*}(m(\gamma(t)))(X_t,Y_t)=0.\eqno{(5.4)}$$
Multiplying  $l_{\scriptscriptstyle (m(\gamma(t)))^{-1}}$ with
both side of (5.4) then  differentiating it at $t=0,$ we get
$(L_{(m_*V)}\pi_{G^*})(X^l_t,Y^l_t)=0.$ Use (2.2) and the fact
that $\{X^l_t.Y^l_t\}=[X_t,Y_t]^l,$ we have
$$\la  m_*V, [X_t,Y_t]^l \ra=
\la \pi_{G^*}^{\sh}(X^l_t),d(i( m_*V)Y^l_t)\ra
 -\la \pi_{G^*}^{\sh}(Y^l_t),d(i( m_*V)X^l_t)\ra=0.$$
Varying $V\in T_p P,$ we have $[X,Y]^l\in {\rm Im}(m_*).$ Use
Proposition 5.1, we immediately have $[\fg_u,\fg_u]\subset \fg_p.$
\qed

\bli{Theorem 5.1} Let $G$ be a compact Poisson Lie group and $P$
an almost homogeneous Poisson $G$-manifold with an equivariant
momentum map ${\rm m}: P\lw G^{*}.$ Then the open $G$-orbit in $P$
is a symplectic fibre bundle over a dressing orbit. In particular,
a compact homogeneous Poisson manifold with an equivariant
momentum map is a symplectic torus bundle over a dressing orbit.
\eli
   \pf  The first part of this theorem follows from the
   discussions at the beginning of this section. If $P$
  is  a compact homogeneous Poisson manifold with an equivariant momentum
map $m,$  let $P={\cal O}_{p}=G\cdot p=G/G_p$ for some $p\in P$
and ${\cal O}_{u}=G/G_u\subset G^*,$ here $ u=m(p).$
 Then $P$ is a symplectic homogeneous fibre bundle over the
dressing orbit ${\cal O}_{u}.$ Since $P$ is compact,  the fibre is
 isomorphic to $G_u/G_p$ and is compact as well.
 By Proposition 5.2, $P$ is a torus bundle over ${\cal
 O}_{u}.$\qed

   At the end of this section we give an example
    of non-transitive almost homogeneous
   Poisson action.
\bli{Example 5.1} {\rm Consider the Lie group
$$SL(2,\R)={\Bigg\{} \left(\begin{array}{cc}
     a & b\\
     c & d
    \end {array}\right):a,b,c,d \in \R,ad-bc=1{\Bigg\}}. $$
We  denote
$$e_{1}=\frac{1}{2}\left(\begin{array}{cc}
      1 & 0\\
      0 & -1
     \end {array}\right), \quad
e_{2}=\frac{1}{2}\left(\begin{array}{cc}
      0 & 1\\
     -1 & 0
     \end{array}\right), \quad
e_{3}=\frac{1}{2}\left(\begin{array}{cc}
      0 & 1\\
      1 & 0
     \end {array}\right)
$$
 a basis of its Lie algebra $sl(2,\R).$    Then
$$[e_{1},e_{2}]=e_{3}, \quad [e_{2},e_{3}]=e_{1}, \quad[e_{3},e_{1}]=-e_{2}.$$
Let $\Lambda= {\lambda}_{1} e_{1} \w e_{2} + {\lambda}_{2} e_{2}
\w e_{3} + {\lambda}_{3} e_{3} \w e_{1}, $ where $
{\lambda}_{1},{\lambda}_{2}, {\lambda}_{3} \in \R.$ Clearly it is
a {\bf r}-matrix of $sl(2,\R).$  Let $\pi (g)
=L_{g*}\Lambda-R_{g*}\Lambda$ with $g\in SL(2, \R).$
  Then $(SL(2, \R),\pi)$ is a  Poisson Lie group. The dual Lie bracket are defined by
   $[\xi,\eta]_{*}=ad^{*}_{{\Lambda} \xi} \eta -
ad^{*}_{{\Lambda} \eta } \xi$   for any  $\xi,  \eta \in
sl^{*}(2,\R).$   Let  $e^{*}_{1},e^{*}_{2},e^{*}_{3}$ be the basis
of $sl^{*}(2,\R)$    dual to $e_{1},e_{2},e_{3}.$
  It is easy to check the Lie brackets on this  basis are given by $$
[e^{*}_{1},e^{*}_{2}]_{*}=-\lambda_{3}e^{*}_{1}-{\lambda}_{2}e^{*}_{2},
\quad
[e^{*}_{2},e^{*}_{3}]_{*}={\lambda}_{1}e^{*}_{2}-{\lambda}_{3}e^{*}_{3},
\quad
[e^{*}_{3},e^{*}_{1}]_{*}={\lambda}_{2}e^{*}_{3}-{\lambda}_{1}e^{*}_{1}.$$
Let $\bpi=h(x_{1},x_{2}){\partial_{x_1}} \w {\partial_{x_1}}.$
Then $({\R}^{2},\bpi)$ is a Poisson manifold for any
$h(x_{1},x_{2})
 \in C^{\infty}({\R}^{2}).$
The natural action of $SL(2,\R)$ on $\R^{2}$ is not transitive and
has two orbits: a fixed point orbit ${\cal O}_{1}=\{0\}$  and an
open orbit ${\cal O}_{2}=\R^{2}-\{0\}.$ In order that the natural
action $\sigma: SL(2, \R)\times ({\R}^{2}, \bpi) \lw
({\R}^{2},\bpi)$ is a Poisson action, we must have
$$ \bpi(gx)={{\sigma}}_{{g}_{*}} \bpi(x)+ {{{\sigma}}^{x}_{*}}\pi(g).
\eqno{(5.5)}
$$
Let $$g= \left(\begin{array}{cc}
      a_{1} & a_{2}\\
      a_{3} &  a_{4}
     \end{array}\right) \in SL(2,\R), \quad\quad
x=\left(\begin{array}{c}
       x_{1}\\
       x_{2}
     \end{array}\right) \in \R^{2},$$
Then (5.5) can be written as
$$
\label{exa2}
\begin{array}{rcl}
&&h(a_{1}x_{1}+a_{2}x_{2}, a_{3}x_{1}+a_{4}x_{2})-h(x_{1},x_{2})\\
&=&\frac{1}{4}\{[({\lambda}_{1}+{\lambda}_{3})a_{1}^{2}-
({\lambda}_{1}-{\lambda}_{3})a_{3}^{2}-2{\lambda}_{2}a_{1}a_{3}-({\lambda}_{1}+{\lambda}_{3})]x_{1}^{2}\\
&&+
[({\lambda}_{1}+{\lambda}_{3})a_{2}^{2}-
({\lambda}_{1}-{\lambda}_{3})a_{4}^{2}-2{\lambda}_{2}a_{2}a_{4}+({\lambda}_{1}-{\lambda}_{3})]x_{2}^{2}\\
&&+[2({\lambda}_{1}+{\lambda}_{3})a_{1}a_{2}
-2({\lambda}_{1}-{\lambda}_{3})a_{3}a_{4}-4{\lambda}_{2}a_{2}a_{3}]x_{1}x_{2}\}\\
\end{array}\eqno{(5.6)}
$$
The solution of (5.6) is
$$
 h(x_{1},x_{2})=\frac{1}{4}({\lambda}_{1}+{\lambda}_{3})x_{1}^{2}-
\frac{1}{4}({\lambda}_{1}-{\lambda}_{3})x_{2}^{2}-\frac{1}{2}{\lambda}_{2}x_{1}x_{2}+c, \quad
\quad c\in \R.
$$
The infinitesimal action of the natural action of $SL(2,\R)$ on $\R^2$ is
$$(e_{1})_{P}=\frac{1}{2} (x_{1}{\partial_{x_{1}}} -x_{2} {\partial_{x_{2}}}),\quad
(e_{2})_{P}=\frac{1}{2} (x_{2}{\partial_{x_{1}}} -x_{1}
{\partial_{x_{2}}} ), \quad (e_{3})_{P}=\frac{1}{2}
(x_{2}{\partial_{x_{1}}} +x_{1} {\partial_{x_{2}}} ).$$ So this
action is a tangential Poisson action if and only if
${\lambda}_{1}, {\lambda}_{2}, {\lambda}_{3}, c$ satisfy
$${{\lambda}_{1}}+{\lambda}_{3}>0,\quad {{\lambda}_{1}}^{2}+
{{\lambda}_{2}}^{2}-{{\lambda}_{3}}^{2}<0, \quad
 c\geq 0.\eqno{(5.7)}$$ For example, if
${\lambda}_{1}={\lambda}_{2}=0$ and ${\lambda}_{3}=4,$ then
$h(x_{1},x_{2})=x_{1}^{2}+x_{2}^{2}+c.$ The natural action $SL(2,
\R)\times ({\R}^{2}, \bpi) \lw ({\R}^{2},\bpi)$ is a tangential
Poisson action if $c\geq 0.$

In the following, we consider the case
 ${\lambda}_{1}={\lambda}_{3}=0$ and ${\lambda}_{2}=2, $  and take
  $h(x_{1},x_{2})=c-x_{1}x_{2}.$
  The dual Poisson Lie group $SL(2,\R)^{*}$ is realized
explicitly  by [Lu]:
$$SL^{*}(2,\R)\stackrel{\rm Def}=SB(2,\R)\stackrel{\rm Def}=
{\Bigg\{} \left(\begin{array}{cc}
     a & b+ic\\
     0 &  a^{-1}
    \end {array}\right):a>0, b,c \in \R,{\Bigg\}}. $$
The dual basis is realized by
$$
e^{*}_{1}=\left(\begin{array}{cc}
      -1 & 0\\
       0 & 1
      \end {array}\right), \quad
e^{*}_{2}=\left(\begin{array}{cc}
      0 & 1\\
      0 & 0
     \end{array}\right), \quad
e^{*}_{3}=\left(\begin{array}{cc}
      0 & i\\
      0 & 0
     \end {array}\right).$$
And the dual  Lie brackets are given by
$$[e^{*}_{1},e^{*}_{2}]_{*}=-2e^{*}_{2}, \quad
[e^{*}_{2},e^{*}_{3}]_{*}=0,       \quad
[e^{*}_{3},e^{*}_{1}]_{*}=2e^{*}_{3}.$$ Noticing that
$(e_{1})^{\circ}\cong {\rm Span}_{\R}\{e_{2}, e_{3}\} $ is an
ideal of $sb(2, \R),$ so $H={\rm Diag}(a,a^{-1})$ is a Poisson Lie
subgroup of $SL(2, \R)$ (\cf. [Y]). Each orbit of the restricted
action of $H$ on ${\R}^{2}$ if of form $\{(x,y)\in
{\R}^{2}:xy=c\},$ where $c$ is any constant. Clearly this action
is Poisson structure preserved and tangential, it has a family of
momentum maps:
$${\rm m}_{H}: {\R}^{2} \lw \R,
(x_{1}, x_{2})\lo \alpha (|c-x_{1}x_{2}|)^{-\frac{1}{2}}, \quad
\alpha\in \R.$$ Clearly they are equivariant momentum maps.} \eli

\section*{Acknowledgement}

 Part of this work was done while I was
 visiting as a guest fellow at
 the Institut f\"ur Mathematik, Ruhr Universit\"at Bochum, Germany.
I would like to thank Prof. A. Huckleberry  and P. Heinzner for
their financial supports. I also wish to thank Prof. Z. J. Liu for
many skillful
 comments and helpful discussions.
\begin{center}
{\bf\Large References}
\end{center}

{\parindent=0pt
\def\toto#1#2{\centerline{\makebox[1.5cm][r] {#1\hss}
\parbox[t]{13cm}{#2}}\vspace{\baselineskip}}

\toto{[AM]}{{\rm  Abraham, R., Marsden, J. E.,} { Foundation of mechanics,
 2nd edition, Benjamin Cummings Reading,} 1978.}

 \toto{[BD]}{{\rm Belavin, A.A., \& Drinfeld, V.G.}
{Triangle equations and simple Lie algebra,}
 {\it Math. Phys. Rev.,} {\bf 4} (1984), 93-165.}

 \toto{[BKM]}{{\rm  Brihaye, Y., Kowalczyk, E. \& Maslanka, P.} {
 Poisson Lie structure on Galelei group,} math. QA0006167.}

\toto{[Dr1]}{{\rm Drinfel¡¯d, V. G.,} { Quntum groups,}
{\it Proc. Cong. Math.,} {\bf Vol1}, Berkeley, 1986.}

\toto{[Dr2]}{ {\rm V. G. Drinfel'd,}
{ On Poisson homogeneous spaces of
 Poisson-Lie groups,}
 {\it Theoret. and Math Phys,} {\bf 95} (1993), 524-525.}

\toto{[Gi]}{{\rm V. Ginzburg,} { Momentum mappings and Poisson
cohomology,} {\it Int. J. Math.,} {\bf 6} (1987), 330-358.}

\toto{[HW]}{ {\rm  Huckleberry, A.\& Wurzbacher, T.,}
{Multiplicity-free complex manifolds,}
 {\it Math. Ann.,} {\bf 286} (1990), 261-280.}

\toto{[HO]}{ {\rm  Huckleberry, A.\& Oeljeklaus, E.,}
{Classification theorems for almost homogeneous spaces,}
 {\it Institute Elie Cartan,} {\bf 9} 1980.}

\toto{[LWX]}{{\rm  Liu, Z.-J., Weinstein, A. \&  X, P.,} { Manin
triples for Lie bialgebroids,} {\it J. Diff. Geom.,} {\bf 45}
(1997), 547-574.}

\toto{[Lu]}{{\rm Lu, J.-H.,}
 {Momentum mappings and reduction of Poisson action,}
S{\'e}minaire Sud-Rhodanien de G{\'e}om{\'e}trie 1989, P. Dazord et
A. Weinstein {\'e}d. MSRI publications, Springer-Verlag, 1990.}

  \toto{[LS]}{{\rm  Levendorskii  \&   Soibelman Y.}
  { Algebras of   functions on compact quantum groups, Schubert cells and quantum tori.}
   {\it Comm. Math. phys.}  {\bf 139} (1991), 141-170.}

\toto{[LW]}{{\rm Lu, J.-H. \& Weinstein, A.,}
 { Poisson Lie groups, dressing
 transformations, and the Bruhat decomposition,} {\it J. Diff. Geom.,}
 {\bf 31} (1990), 501-526.}

\toto{[STS]}{{\rm  Semenov-Tian-Shansky, M. A.,}
 { Dressing transformations  and Poisson group actions,}
 {\it Publ. Res. Inst. Sci. Kyoto Univ.}
 {\bf 21} (1985), 1237-1260.}

     \toto{[V]}{{\rm  I. Vaisman,}
 {Lectures on the geometry of Poisson manifolds,}
 {\bf 118}, Birkhauser, Basel, 1994.}

  \toto{[W]}{{\rm A. Weinstein,} {Some remarks on dressing transformations,}
 {\it J. Math. Fac. Sci. Univ. Tokyo.,}  {\bf 36} (1988), 163-167.}

  \toto{[Y]}{{\rm Q.-L. Yang,} {Tangential Poisson actions and Poisson reductions,}
 {\it Adv. Math. Chin..} {\bf 31} (2002), 127-134.}

\toto{[Z]}{{\rm  Zakrzewski, S.,}
 {Poisson  structures on Poincar\'e group,}
{\it Comm. Math. Phys.} {\bf 185} (1997), 285-311.}

}

\end{document}